\documentclass[12pt]{amsart}

\usepackage{amsthm,latexsym,amssymb,amsfonts,amsmath,url,mathtools,mathrsfs}
	
\usepackage{placeins}
\usepackage{xcolor}
\usepackage[vcentermath,enableskew]{youngtab}
\usepackage{color}
\usepackage[noadjust]{cite}
\usepackage[margin=1in]{geometry}
\usepackage{enumitem, bm}
\usepackage[unicode,psdextra,linkcolor=purple,colorlinks=true,citecolor=green]{hyperref}
\usepackage{tikz-cd}

\usepackage{ytableau}

\usepackage{graphicx}
\graphicspath{{./images/}}

\usepackage{tikz}
\usetikzlibrary{decorations.markings}

\usetikzlibrary{calc, shapes, backgrounds,arrows,positioning,plotmarks}
\tikzset{>=stealth',
	head/.style = {fill = white, text=black},
	plaque/.style = {draw, rectangle, minimum size = 10mm}, 
	pil/.style={->,thick},
	junct/.style = {draw,circle,inner sep=0.5pt,outer sep=0pt, fill=black}
}

\newtheorem{theorem}{Theorem}[section]
\def\sectionautorefname~#1\null{Section #1\null}
\def\subsectionautorefname~#1\null{Subsection #1\null}
\newtheorem{proposition}[theorem]{Proposition}

\newtheorem{lemma}[theorem]{Lemma}

\newtheorem{corollary}[theorem]{Corollary}

\theoremstyle{definition}
\newtheorem{definition}[theorem]{Definition}

\newtheoremstyle{named}{}{}{\itshape}{}{\bfseries}{.}{.5em}{#1\thmnote{ #3}}
\theoremstyle{named}

\numberwithin{equation}{section}

\newcommand{\LIS}{\mathrm{LIS}}
\newcommand{\Ind}{\mathrm{Ind}}
\newcommand{\Res}{\mathrm{Res}}
\newcommand{\SYT}{\mathrm{SYT}}

\newcommand{\RR}{\mathbb R}
\newcommand{\RRR}{\mathcal R}
\newcommand{\PPP}{\mathcal P}
\newcommand{\CCC}{\mathcal C}
\newcommand{\JJJ}{\mathcal J}
\newcommand{\Prob}{\operatorname{Prob}}
\newtheorem*{thm:Chen}{Theorem \ref{thm:Chen}}
\title{log-concavity of longest increasing subsequence}

\author{Yupeng Li}
\address{Department of Mathematics, Michigan State University, East Lansing, MI 48824, USA}
\email{yupengli@msu.edu}

\date{\today}

\subjclass[2020]{Primary: 05A20, 60J10; Secondary: 05A05, 05E10}

\begin{document}

	\begin{abstract}
	    We prove a conjecture proposed by Chen \cite{Chen08} that, for every fixed $n$, the sequence $(a_{n,k})_{k=1}^n$ is log-concave, where $a_{n,k}$ counts permutations of $[n]$ whose longest increasing subsequence has length $k$.
	\end{abstract}

	\maketitle	

\section{Introduction}
For $\pi\in S_n$, write $\LIS(\pi)$ for the length of a longest strictly increasing subsequence of $\pi$, and set 
\[
a_{n,k} = |\{\pi\in S_n : \LIS(\pi) = k\}|, \quad 1\leq k \leq n.
\]
Chen conjectured that this sequence is log-concave for every fixed $n$ \cite[Conjecture 1.1]{Chen08}. B\'ona, Lackner and Sagan \cite{BLS17} studied this conjecture, its involution analogue and restrictions to specified tableau shapes. Baik, Deift, and Johansson \cite{BDJ99} proved that, with suitable scaling, the sequence converges to the Tracy-Widom distribution as $n \rightarrow \infty$. Rhoades \cite{BR24} constructed a graded quotient ring whose Hilbert function is the reversed sequence of $a_{n,k}$. Novak and Rhoades \cite{NR24} proposed a representation-theoretic strengthening involving Kronecker coefficients. More recently, Baslingker, Krishnapur and Madiman \cite{BKM2026} showed that for every fixed $j$, Chen's inequalities hold for $k = n-j,\dots, n$ once $n$ is sufficiently large. Xia \cite{Xia2026} showed that Chen's conjecture is true when $n \geq 10^{10^{65}}$. For further background, see \cite{Romik15}.

Here we prove Chen's conjecture in full generality and show that the inequalities are strict. The idea of the proof follows closely to the method H\"aggstr\"om \cite{HÄGGSTRÖM} used to establish stochastic domination between the conditional distributions of a vector of independent Bernoulli random variables given successive values of their sum. 

The key intermediate result is \autoref{thm:couplingRnk}, which gives a coupling of Plancherel measures conditioned on first-row lengths $k$ and $k+1$, supported on pairs of partitions $(\lambda,\nu)$ satisfying $\nu_i\le\lambda_i+1$ for every $i$. For such pairs, we compare the probabilities that the first row grows or shrinks under the down-up Markov chain. The coupling transfers these pointwise comparisons to inequalities between the corresponding conditional expectations. Combining these inequalities with the reversibility of the chain yields strict log-concavity.
\begin{theorem}[{\cite[Conjecture 1.1]{Chen08}}]\label{thm:Chen}
    For every integer $n\geq 3$ and every $1< k < n$,
    \[
    a_{n,k}^2 > a_{n,k-1}a_{n,k+1}.
    \]
    In particular, the sequence $(a_{n,1},\dots,a_{n,n})$ is strictly log-concave.
\end{theorem}

\bigskip

\noindent
\textbf{Outline.}
The basic outline of this paper is as follows. In \autoref{sec:probability}, we recall basic terminology and results from probability and finite Markov chains. In \autoref{sec:combinatorics}, we introduce notation for partitions and standard Young tableaux, review the branching identities, and define the Plancherel measure. In \autoref{sec:coupling}, we construct couplings of Plancherel measures conditioned on adjacent first-row lengths, using Markov chains that preserve a specified relation between partitions. Finally, in \autoref{sec:logconcavity}, we use these couplings to compare the probabilities that the first row grows or shrinks. Combining these comparisons with the reversibility of the Markov chain yields the strict log-concavity of $(a_{n,1},\ldots,a_{n,n})$.

\section{Background and Notation: Probability and finite Markov chains}\label{sec:probability}

This section recalls the probability terminology used in the proofs. All state spaces are finite, and all Markov chains evolve in discrete time. The standard reference is Levin, Peres, and Wilmer~\cite{LPW09}.

\subsection{Distributions, expectations, and conditioning}

Let $S$ be a finite nonempty set. A \emph{probability distribution} on $S$ is a function $p:S\to[0,1]$ such that
\[
 \sum_{x\in S}p(x)=1.
\]

Let $X$ have distribution $p$. For a function $g:S\to\RR$, its \emph{expectation} under $p$ is the weighted sum
\[
 \mathbb E_p[g]=\mathbb E[g(X)]=\sum_{x\in S}p(x)g(x).
\]

For $A\subset S$, set $p(A) = \sum_{x\in A}p(x)$. If $p(A)>0$, the distribution \emph{conditioned on} $A$ assigns probability $p(x)/p(A)$ to each $x\in A$ and zero to $x\notin A$. For example, the measure $\pi_{n,k}$ defined below is obtained by conditioning Plancherel measure on the first-row length being $k$.

\subsection{Joint distributions and couplings}

Let $p$ and $q$ be probability distributions on finite sets $S$ and $T$. A \emph{joint distribution} of random variables $X$ and $Y$ is a probability distribution $\Gamma$ on $S\times T$, where
$\Gamma(x,y)=\Prob(X=x,Y=y)$. Its two \emph{marginal distributions} are the individual distributions of $X$ and $Y$, obtained by summing over the other variable. A \emph{coupling} of $p$ and $q$ is a joint distribution satisfying
\[
 \sum_{y\in T}\Gamma(x,y)=p(x),\qquad
 \sum_{x\in S}\Gamma(x,y)=q(y).
\]
The product $\Gamma(x,y)=p(x)q(y)$ gives the independent coupling; other couplings allow dependence between the two choices. The \emph{support} of $\Gamma$ is the set of pairs with positive probability. A coupling is \emph{supported on} a relation $F\subseteq S\times T$ if $\Gamma(x,y)=0$ outside $F$. Equivalently, $\Prob((X,Y)\in F)=1$.

Couplings also compare expectations. For functions $g:S\to\RR$ and $h:T\to\RR$, if $g(x)\ge h(y)$ whenever $\Gamma(x,y)>0$, then the coupling identities give
\[
 \mathbb E_p[g]-\mathbb E_q[h]
 =\sum_{x\in S}\sum_{y\in T}\Gamma(x,y)\bigl(g(x)-h(y)\bigr)
 \ge0.
\]
The inequality is strict if $g(x)>h(y)$ for at least one pair with positive probability.

\subsection{Transition matrices and Markov chains}

A \emph{transition matrix}, also called a \emph{Markov kernel}, on $S$ is a function $H:S\times S\to[0,1]$ with
\[
 \sum_{y\in S}H(x,y)=1, \qquad x\in S.
\]
A Markov chain with transition matrix $H$ is a sequence $X_0,X_1,\ldots$ whose next state is chosen according to the row $H(X_t,\cdot)$. More precisely, for every history with positive
probability,
\[
 \Prob(X_{t+1}=y\mid X_0=x_0,\ldots,X_t=x)=H(x,y).
\]
Thus the distribution of the next state depends only on the current state, and the same transition rule is used at every time. If the current distribution is $p$, the distribution after one step is the
row vector $pH$, with
\[
 (pH)(y)=\sum_{x\in S}p(x)H(x,y).
\]

\subsection{Stationarity, reversibility, and irreducibility}

A distribution $p$ is \emph{stationary} for $H$ if $pH=p$, or equivalently
\[
 p(y)=\sum_{x\in S}p(x)H(x,y),\qquad y\in S.
\]
Starting the chain with distribution $p$ then gives the same distribution at every time.

The matrix $H$ is \emph{reversible with respect to} a distribution $p$ if it satisfies the \emph{detailed-balance equations}
\[
 p(x)H(x,y)=p(y)H(y,x),\qquad x,y\in S.
\]
These equations equate the probability of each transition with that of its reverse when the initial state has distribution $p$. They imply stationarity, since
\[
 \sum_{x\in S}p(x)H(x,y)
 =p(y)\sum_{x\in S}H(y,x)=p(y).
\]

The chain is \emph{irreducible} if any state can be reached from any other by a path of transitions with positive probability. Equivalently, for every $x,y\in S$, some integer $t\ge0$ satisfies $H^t(x,y)>0$. Every transition matrix on a finite nonempty set has a stationary distribution. If it is irreducible, that distribution is unique and assigns positive probability to every state. 

\subsection{Coupled transitions}\label{subsec:coupledtrans}

Let $H_1$ and $H_2$ be transition matrices on $S$ and $T$. A transition matrix $G$ on $S\times T$ \emph{couples} their transitions if, for every current pair $(x,y)$ and all $x'\in S$, $y'\in T$,
\[
 \begin{aligned}
 \sum_{y'\in T}G((x,y),(x',y'))&=H_1(x,x'),\\
 \sum_{x'\in S}G((x,y),(x',y'))&=H_2(y,y').
 \end{aligned}
\]
The first identity specifies the distribution of the first next state by summing over all possibilities for the second next state; the second identity has the analogous meaning. The two chains therefore follow the transition rules $H_1$ and $H_2$, respectively, while their moves may be dependent.

The transitions \emph{preserve} a nonempty relation $F\subseteq S\times T$ if $G((x,y),(x',y'))=0$ whenever $(x,y)\in F$ and $(x',y')\notin F$. The restriction of $G$ to $F$ is then the transition matrix of a finite Markov chain and has a stationary distribution. Equivalently, the transition matrix can be defined directly on $F$, with the marginal sums restricted to next pairs belonging to $F$. Its stationary distribution has marginals that are stationary for $H_1$ and $H_2$, by the two identities above. When $H_1$ and $H_2$ are irreducible, uniqueness identifies these marginals as their respective stationary distributions.

This method of constructing a coupling through a stationary distribution has a classical precedent in Holley's proof of \cite[Theorem~6]{Holley74}. Holley constructed a coupled continuous-time Markov chain whose transitions preserve the order on a finite distributive lattice. Its limiting stationary distribution gives a coupling of the prescribed distributions supported on ordered pairs. Leskel\"a \cite{Leskela10} developed a general theory of couplings of Markov processes that preserve a specified relation and gives conditions for their stationary distributions to admit a coupling supported on that relation.

\section{Background and Notation: Combinatorics}\label{sec:combinatorics}
Write $f^\lambda = |\SYT(\lambda)|$ for the number of standard Young tableaux of shape $\lambda$ in English notation. Let $\PPP_n$ denote the set of partitions of $n$, with $\PPP_0 = \{\varnothing\}$, and let $\PPP_{n,k}$ denote partitions of $n$ whose first row has length $k$. We write $\nu\nearrow \lambda$ when $\lambda$ is obtained from $\nu$ by adding one box and write $\nu\searrow \lambda$ when $\lambda$ is obtained from $\nu$ by deleting one box. We extend partitions by zeros beyond their nonzero parts, and write $e_i$ for the $i$-th unit vector. Therefore, if $\nu\nearrow \lambda$, then $\lambda = \nu +e_i$ for some $i$. Similarly, if $\nu\searrow \lambda$, then $\lambda = \nu - e_i$ for some $i$.

The standard branching identities are 
\begin{align}
    f^\lambda &= \sum_{\nu\nearrow \lambda}f^\nu,\qquad \lambda\in \PPP_n, n\geq 1\label{eq:godown},\\
    (n+1)f^\lambda &= \sum_{\lambda\nearrow\eta}f^\eta,\qquad \lambda\in \PPP_n, n\geq 0\label{eq:goup}.
\end{align}
We set $f^\varnothing = 1$. Let $V^\lambda$ represent the complex Specht module. Then \autoref{eq:godown} follows from
\[
\Res_{S_{n-1}}^{S_n}V^\lambda \cong \bigoplus_{\nu\nearrow\lambda}V^\nu,
\] and \autoref{eq:goup} follows from 
\[\Ind_{S_n}^{S_{n+1}}V^\lambda \cong \bigoplus_{\lambda\nearrow\eta}V^\eta.\]
Robinson–Schensted gives the identities
\[
a_{n,k} = \sum_{\lambda\in \PPP_{n,k}}(f^\lambda)^2,\text{ and } n! = \sum_{\lambda\in \PPP_n}(f^\lambda)^2.
\]
Let $\pi_n(\lambda) = \frac{(f^\lambda)^2}{n!}$ be the \textit{Plancherel measure} \cite{BOO00} on $\PPP_n$.
Conditioning Plancherel measure on the first-row length gives
\[
\pi_{n,k}(\lambda) = \frac{(f^\lambda)^2}{a_{n,k}},\qquad \lambda\in \PPP_{n,k}.
\]
\section{Relation preserved by random walk}\label{sec:coupling}
For $\lambda \in \PPP_n$, we define removal and addition rows probabilities by
\begin{align}
    D_i(\lambda) &= \frac{f^{\lambda-e_i}}{f^\lambda},\qquad n\geq 1,\\
    U_i(\lambda) &= \frac{f^{\lambda+e_i}}{(n+1)f^\lambda},\qquad n\geq 0.
\end{align}
We assign probability zero to illegal moves. The branching identities
\autoref{eq:godown} and \autoref{eq:goup} give
\[
    \sum_i D_i(\lambda)=1,
    \qquad
    \sum_i U_i(\lambda)=1. 
\]
Thus $D_i(\lambda)$ and $U_i(\lambda)$ are probability distributions on removal rows and addition rows, respectively.
\subsection{Random walk on \texorpdfstring{$\PPP_n$}{P\_n}} For fixed $n\geq 1$, let $K_n$ be the transition matrix of a Markov chain on $\PPP_n$ that removes a box according to $D$ and then adds a box according to $U$ from the resulting partition. Write $K_n(\lambda,\eta)$ for the transition probability from $\lambda$ to $\eta$. These transitions and their reversibility are discussed in \cite[Proposition 2.1(4)]{Fulman05}; the formulas needed here are derived explicitly below.
\begin{lemma}\label{lem:reversibility}
    Let $c(\lambda,\eta) = |\{\rho \in \PPP_{n-1} \,:\, \lambda\searrow \rho\nearrow \eta\}|$. Then $K_n(\lambda,\eta) = c(\lambda,\eta)\frac{f^\eta}{nf^\lambda}$ and
    \begin{equation}\label{eq:reversibility}
        \pi_n(\lambda)K_n(\lambda, \eta) = \pi_n(\eta)K_n(\eta,\lambda) = \frac{c(\lambda,\eta)f^\lambda f^\eta}{n n!}.
    \end{equation}
    Consequently $K_n$ is reversible with respect to $\pi_n$, and $\pi_n$ is a stationary distribution for $K_n$.
\end{lemma}
\begin{proof}
    Each path $\lambda\searrow\rho\nearrow\eta$ has probability
    \[
    \frac{f^\rho}{f^\lambda}\cdot \frac{f^\eta}{nf^\rho} = \frac{f^\eta}{nf^\lambda}.
    \]
    Summing over all paths gives $K_n(\lambda,\eta) = c(\lambda,\eta)\frac{f^\eta}{nf^\lambda}$. Also, we have $c(\lambda,\eta) = c(\eta,\lambda)$ and $\pi_n(\lambda) = \frac{(f^\lambda)^2}{n!}$ by definition. Therefore, 
    \[
     \pi_n(\lambda)K_n(\lambda, \eta) =
     \frac{c(\lambda,\eta)f^\lambda f^\eta}{n n!} 
     =\frac{c(\eta,\lambda)f^\lambda f^\eta}{n n!} =
     \pi_n(\eta)K_n(\eta,\lambda).
    \]
\end{proof}
\begin{lemma}\label{lem:explicitformula}
    Let $\lambda \in \PPP_n$, choose $L\geq n+1$, append zeros to $\lambda$ so that it has $L$ rows, and set $b_i = \lambda_i + L -i$. For $1\leq i\leq L$,
\begin{align}
    D_i(\lambda) &= \frac{b_i}{n}\prod_{j\neq i}\frac{b_i-b_j-1}{b_i-b_j},\quad n\geq 1,\\
    U_i(\lambda) &= \frac{1}{b_i+1}\prod_{j\neq i}\frac{b_i-b_j+1}{b_i-b_j},\quad n\geq 0.
\end{align}
The formulas vanish at illegal moves.
\end{lemma}
\begin{proof}
    The Vandermonde form of the hook-length formula is
    \[
    f^\lambda = n! \frac{\prod_{1\leq k < j\leq L}(b_k-b_j)}{\prod_{k=1}^L b_k!}.
    \]
    For a legal removal or addition in row $i$, replace $b_i$ by $b_i-1$ or $b_i+1$, respectively. This gives
    \[
    D_i(\lambda) = \frac{(n-1)!b_i!}{n! (b_i-1)!}\prod_{i < j}\frac{b_i - 1 -b_j}{b_i -b_j}\prod_{k < i}\frac{b_k - b_i + 1}{b_k - b_i} = \frac{b_i}{n}\prod_{j\neq i}\frac{b_i-b_j-1}{b_i-b_j},
    \]
    and a similar calculation gives the expression for $U_i(\lambda)$.

    Since $(b_i)$ is a strictly decreasing nonnegative sequence, it suffices to show that the numerators of $D_i(\lambda)$ and $U_i(\lambda)$ are zero for illegal moves. For an illegal removal in row $i < L$, $\lambda_i = \lambda_{i+1}$, so $b_i - b_{i+1} - 1= 0$ which implies $D_i(\lambda) = 0$. An illegal removal in the last appended row implies $b_L = 0$ so the numerator vanishes and $D_L(\lambda) = 0$. For an illegal addition in row $i > 1$, $\lambda_i = \lambda_{i-1}$, so $b_i - b_{i-1} + 1 = 0$ which implies $U_i(\lambda) = 0$. 
\end{proof}
\subsection{Relation on \texorpdfstring{$\PPP_n$}{P\_n}} We introduce a relation on $\PPP_n$ and construct couplings of distributions on $\PPP_n$ supported on this relation.
\begin{definition}
    For $\lambda,\nu\in \PPP_n$, we write $\lambda\RRR \nu$ if $\nu_i\leq \lambda_i + 1$ for all $i$. We say row $i$ is tight if $\nu_i = \lambda_i + 1$.
\end{definition}
We use $\RRR$ for the relation; transitivity is neither asserted nor needed. In particular, the results of
this section apply to every related pair of the same size.
\begin{lemma}\label{lem:DU}
    Suppose $\lambda, \nu \in \PPP_n$ and $\lambda\RRR\nu$. At every tight row $i$, \[
    D_i(\lambda) \leq D_i(\nu),\qquad U_i(\lambda) \geq U_i(\nu).
    \]
    If row one is tight, then $U_1(\lambda) > U_1(\nu)$.
\end{lemma}
\begin{proof}
    Fix $L \geq n+1$, and write $b_j = \lambda_j + L - j$ and $c_j = \nu_j + L - j$. At a tight row $i$, we have $c_i = b_i + 1$. Since $c_j \leq b_j + 1$ for every $j$, it follows that $c_i-c_j \geq b_i-b_j$. Both $(b_j)$ and $(c_j)$ are strictly decreasing sequences, so
    \[
    0 > c_i-c_j \geq b_i-b_j
    \qquad\text{for }j<i,
    \]
    and
    \[
    c_i-c_j \geq b_i-b_j > 0
    \qquad\text{for }j>i.
    \]
    Since $\lambda_L=\nu_L=0$, a tight row satisfies $i<L$. Suppose first that $b_i-b_{i+1}\geq 2$. By \autoref{lem:explicitformula},
    \[
    \frac{D_i(\nu)}{D_i(\lambda)}
    =
    \frac{c_i}{b_i}\prod_{j<i}\frac{(c_i-c_j-1)(b_i-b_j)}{(b_i-b_j-1)(c_i-c_j)}\prod_{j>i}\frac{(c_i-c_j-1)(b_i-b_j)}{(b_i-b_j-1)(c_i-c_j)}
    \geq 1.
\]
If $b_i-b_{i+1}=1$, the preceding ratio is undefined
because its denominator vanishes. In this case,
$b_i-b_{i+1}-1=0$ implies $D_i(\lambda)=0$, and hence
$D_i(\lambda)\leq D_i(\nu)$.

A similar argument gives $U_i(\lambda)\geq U_i(\nu)$ when $i\geq 2$. When $i=1$, all factors in $U_1$ are strictly positive, and
\[
\frac{U_1(\lambda)}{U_1(\nu)}
=
\frac{c_1+1}{b_1+1}\prod_{j>1}\frac{(b_1-b_j+1)(c_1-c_j)}{(c_1-c_j+1)(b_1-b_j)}
>1.
\]
\end{proof}
Now let $(p_i)$ and $(q_i)$ be probability distributions on the same finite index set, and choose numbers $0 \leq t_i \leq \min(p_i,q_i)$. Set $s = 1 - \sum_i t_i$ and define
\begin{equation}\label{eq:coupling}
    \mathcal{C}_{ij} = \delta_{ij}t_i + (p_i-t_i)(q_j-t_j)/s,\qquad s>0,
\end{equation}

and $\mathcal{C}_{ij} = \delta_{ij}t_i$ when $s = 0$, where $\delta_{ij}$ is the Kronecker delta. It is not hard to check that $\CCC = (\CCC_{ij})$ is a coupling of $(p_i)$ and $(q_i)$. This is a slight generalization of the coupling in \cite[Section 4.2, Proposition 4.7]{LPW09} where the authors only consider $t_i = \min(p_i,q_i)$.

\begin{proposition}\label{prop:Cmatrix}
    Let $\lambda,\nu\in \PPP_n$ satisfy $\lambda \RRR \nu$. Append zeros so that both partitions have $L\geq n+1$ rows.
    \begin{enumerate}
        \item If $n \geq 1$, the removal distributions $(D_i(\lambda))_{i=1}^L$ and $(D_i(\nu))_{i=1}^L$ admit a coupling $\CCC^{-}_{\lambda,\nu}$ supported on pairs of legal removal rows $(i,j)$ satisfying 
        \[
        (\lambda-e_i)\RRR (\nu-e_j).
        \]
        \item The addition distributions $(U_i(\lambda))_{i=1}^L$ and $(U_i(\nu))_{i=1}^L$ admit a coupling $\CCC^{+}_{\lambda,\nu}$ supported on pairs of legal addition rows $(i,j)$ satisfying 
        \[
        (\lambda+e_i)\RRR (\nu+e_j).
        \]
    \end{enumerate}
\end{proposition}
\begin{proof}
    We need to construct $\CCC^{-}_{\lambda,\nu}$ satisfying 
    \[
    \sum_{j=1}^L (\mathcal{C}^-_{\lambda,\nu})_{ij} = D_i(\lambda), \qquad \sum_{i=1}^L (\mathcal{C}^-_{\lambda,\nu})_{ij} = D_j(\nu),
    \] and
    \[
    (\mathcal{C}^-_{\lambda,\nu})_{ij} > 0 \text{ implies } (\lambda-e_i)\RRR (\nu-e_j).
    \]
    We suppress $(\lambda,\nu)$ in the notation since they are fixed. Let $T$ be the set of tight rows. Since we have the relation $\lambda\RRR\nu$, set $p_i = D_i (\lambda), q_i = D_i(\nu)$, and $t_i = p_i$ when $i\in T$ and $t_i = 0$ otherwise. By \autoref{lem:DU}, $p_i\leq q_i$ on tight rows, so these choices also satisfy $0\leq t_i\leq\min(p_i,q_i)$. Let $\mathcal{C}^{-}$ be the matrix defined by \autoref{eq:coupling}. It has the required row and column sums. Moreover, $\CCC_{ij}^- = 0$ if $i\in T$ and $j\neq i$. Now we assume $\CCC_{ij}^- > 0$. If $i=j$, subtracting one from the same row preserves the relation. If $i\neq j$, then $i\notin T$ and $\nu_i < \lambda_i +1$. Thus, $(\nu-e_j)_k \leq (\lambda-e_i)_k+1$ for all $k$ and hence $(\lambda - e_i) \RRR (\nu -e_j)$.

    For $\mathcal{C}^+$, define $p_i = U_i (\lambda), q_i = U_i(\nu)$, and $t_i = q_i$ when $i\in T$ and $t_i = 0$ otherwise. The rest of the proof is similar.
\end{proof}
\begin{corollary}\label{cor:probdistributeonpairs}
        Let $n \geq 1$ and $\lambda,\nu\in \PPP_n$ satisfy $\lambda \RRR \nu$. The transition distributions $K_n(\lambda,\cdot)$ and $K_n(\nu,\cdot)$ admit a coupling $Q_{\lambda,\nu}$ supported on 
        \[
        \{(\widehat{\lambda},\widehat{\nu})\in \PPP_n\times \PPP_n\,:\,\widehat{\lambda}\RRR\widehat{\nu}\}.
        \]
\end{corollary}
\begin{proof}
Fix $L\geq n+1$, and let all row indices range from $1$ to $L$.
Set
\[
S=\{(i,j):(\CCC^-_{\lambda,\nu})_{ij}>0\},
\qquad
I_\lambda=\{i:D_i(\lambda)>0\}.
\]
By \autoref{prop:Cmatrix}, every $(i,j)\in S$
gives legal, related partitions $\lambda-e_i$ and $\nu-e_j$. Now we define 
\[
Q_{\lambda,\nu}(\widehat{\lambda},\widehat{\nu}) = 
\sum_{(i,j)\in S}(\CCC^-_{\lambda,\nu})_{ij}\bigg(\sum_{\substack{i',j'\\ \lambda-e_i+e_{i'} =\widehat{\lambda}\\ \nu-e_j+e_{j'} = \widehat{\nu}}}(\CCC^+_{\lambda-e_i,\nu-e_j})_{i'j'}\bigg).
\]
The sum above may contain multiple summands when $\widehat{\lambda} =\lambda$ or $\widehat{\nu} = \nu$. To prove the statement, we need to check the following:
\begin{enumerate}
        \item $\displaystyle\sum_{\mu\in \PPP_n}Q_{\lambda,\nu}(\widehat{\lambda},\mu) = K_n(\lambda,\widehat{\lambda})$,
        \item $\displaystyle\sum_{\mu\in \PPP_n}Q_{\lambda,\nu}(\mu,\widehat{\nu}) = K_n(\nu,\widehat{\nu})$, and
        \item If $Q_{\lambda,\nu}(\widehat{\lambda},\widehat{\nu})>0$ then $\widehat{\lambda}\RRR\widehat{\nu}$.
\end{enumerate}
For the first marginal,
\begin{align}
    \sum_{\mu\in \PPP_n}Q_{\lambda,\nu}(\widehat{\lambda},\mu)  
    &= \sum_{\mu\in \PPP_n}\Bigg(\sum_{(i,j)\in S}(\CCC^-_{\lambda,\nu})_{ij}\bigg(\sum_{\substack{i',j'\\ \lambda-e_i+e_{i'} =\widehat{\lambda}\\ \nu-e_j+e_{j'} = \mu}}(\CCC^+_{\lambda-e_i,\nu-e_j})_{i'j'}\bigg)\Bigg)\\
    &= \sum_{\mu\in \PPP_n} \sum_{(i,j)\in S} (\CCC^-_{\lambda,\nu})_{ij}\cdot \bigg(\sum_{\substack{i',j'\\\lambda-e_i+e_{i'} =\widehat{\lambda}\\ \nu-e_j+e_{j'} = \mu}}(\CCC^+_{\lambda-e_i,\nu-e_j})_{i'j'}\bigg)\\
    &= \sum_{(i,j)\in S}(\CCC^-_{\lambda,\nu})_{ij} \cdot\bigg(\sum_{\substack{i'\\ \lambda-e_i+e_{i'} = \widehat{\lambda}}}\sum_{j'=1}^L(\CCC^+_{\lambda-e_i,\nu-e_j})_{i'j'}\bigg)\\
    &= \sum_{(i,j)\in S}(\CCC^-_{\lambda,\nu})_{ij} \cdot\bigg( \sum_{\substack{i'\\ \lambda-e_i+e_{i'} = \widehat{\lambda}}}U_{i'}(\lambda-e_i)\bigg)\label{eq:sumoverj'}\\
    &= \sum_{i\in I_\lambda}\bigg(\sum_{j:(i,j)\in S} (\CCC^{-}_{\lambda,\nu})_{ij}\bigg) \cdot \bigg(\sum_{\substack{i'\\ \lambda-e_i+e_{i'} = \widehat{\lambda}}}U_{i'}(\lambda-e_i)\bigg)\label{eq:movejin}\\
    &= \sum_{i\in I_\lambda}\bigg(D_i(\lambda) \cdot \sum_{\substack{i'\\ \lambda-e_i+e_{i'} = \widehat{\lambda}}}U_{i'}(\lambda-e_i)\bigg)\label{eq:sumoverj}\\
    &= K_n(\lambda,\widehat{\lambda}).
\end{align}
\autoref{eq:sumoverj'} and \autoref{eq:sumoverj} follow from \autoref{prop:Cmatrix}. \autoref{eq:movejin} follows because the second sum is independent of $j$.

The identity $\displaystyle\sum_{\mu\in \PPP_n}Q_{\lambda,\nu}(\mu,\widehat{\nu}) = K_n(\nu,\widehat{\nu})$ follows by the same calculation.

If $Q_{\lambda,\nu}(\widehat{\lambda},\widehat{\nu})>0$, then there exist $i,j,i',j'$ such that $(\CCC^-_{\lambda,\nu})_{ij}$ and $(\CCC^+_{\lambda-e_i,\nu-e_j})_{i'j'}$ are both greater than $0$. By \autoref{prop:Cmatrix}, $(\lambda-e_i)\RRR(\nu-e_j)$. Applying \autoref{prop:Cmatrix} again to the pair $\lambda-e_i,\nu-e_j\in \PPP_{n-1}$ gives $(\lambda-e_i+e_{i'})\RRR(\nu-e_j+e_{j'})$.

Finally, $Q_{\lambda,\nu}\geq 0$ and $\displaystyle \sum_{\widehat{\lambda},\widehat{\nu}}Q_{\lambda,\nu}(\widehat{\lambda},\widehat{\nu}) = \sum_{\widehat{\lambda}}K_n(\lambda,\widehat{\lambda}) = 1$. Thus $Q_{\lambda,\nu}$ is a joint probability distribution.
\end{proof}
\subsection{Random walk on \texorpdfstring{$\PPP_{n,k}$}{Pnk}} Now we restrict the random walk to the set of partitions of $n$ whose first row has length $k$.
\begin{definition}\label{def:restoremove}
    For $n\geq 1$ and $1\leq k \leq n$, define the transition matrix $K_{n,k}$ on $\PPP_{n,k}$ by
    \[
    K_{n,k}(\lambda,\nu)=
    \begin{cases}
      K_n(\lambda,\nu), &\nu\neq \lambda,\\
      K_n(\lambda,\lambda) + \sum_{\mu \in \PPP_n\setminus\PPP_{n,k}}K_n(\lambda,\mu), &\nu=\lambda.
    \end{cases}
    \]
\end{definition}
Thus a $K_{n}$ transition whose outcome lies outside $\PPP_{n,k}$ is replaced by a transition that stays at the current partition.
\begin{lemma}\label{lem:Kmk}
    For $n\geq 1$ and $1\leq k\leq n$, the transition matrix $K_{n,k}$ is irreducible and reversible with respect to $\pi_{n,k}$. Consequently, $\pi_{n,k}$ is its unique stationary distribution.
\end{lemma}
\begin{proof}
    By \autoref{lem:reversibility}, for distinct partitions $\lambda$ and $\nu$, \[(f^\lambda)^2 K_{n,k}(\lambda,\nu) = (f^\nu)^2 K_{n,k}(\nu,\lambda).\] Dividing by $a_{n,k}$ gives \[\pi_{n,k}(\lambda)K_{n,k}(\lambda,\nu) = \pi_{n,k}(\nu)K_{n,k}(\nu,\lambda).\]
    For $\lambda = \nu$, this equality is automatic. To verify stationarity, fix $\nu\in \PPP_{n,k}$ and sum over $\lambda$:
    \[
    \sum_{\lambda\in \PPP_{n,k}}\pi_{n,k}(\lambda)K_{n,k}(\lambda,\nu) = \pi_{n,k}(\nu)\sum_{\lambda\in \PPP_{n,k}} K_{n,k}(\nu,\lambda)= \pi_{n,k}(\nu)\sum_{\lambda\in \PPP_n}K_n(\nu,\lambda) = \pi_{n,k}(\nu).
    \]
    
    For irreducibility, it suffices to show that every $\lambda\in \PPP_{n,k}$ can reach the hook $(k,1^{n-k})$ with positive probability. While there is a row $i\geq 2$ containing more than one box, choose the last such row. Let $\ell$ be the number of nonzero rows of the current partition, and move the last box of row $i$ to a new row $\ell+1$. Each intermediate step $\lambda\rightarrow \lambda-e_i+e_{\ell+1}$ has positive probability given by $D_i(\lambda)U_{\ell+1}(\lambda-e_i)>0$. Every transition used above has a positive reverse transition $\lambda-e_i+e_{\ell+1}\rightarrow \lambda$ with probability $D_{\ell+1}(\lambda-e_i+e_{\ell+1})U_i(\lambda-e_i)$.

    To prove uniqueness, let $\theta$ be any stationary probability distribution for $K_{n,k}$. Since $\pi_{n,k}(\lambda)>0$ for all $\lambda \in \PPP_{n,k}$, define $r(\lambda) = \frac{\theta(\lambda)}{\pi_{n,k}(\lambda)}$. Stationarity of $\theta$ implies
    \begin{align}
    r(\nu) &= \frac{1}{\pi_{n,k}(\nu)}\sum_{\lambda\in \PPP_{n,k}}\theta(\lambda)K_{n,k}(\lambda,\nu)\\
    &= \sum_{\lambda\in\PPP_{n,k}}\frac{r(\lambda)\pi_{n,k}(\lambda)}{\pi_{n,k}(\nu)}K_{n,k}(\lambda,\nu)\\
    &= \sum_{\lambda\in\PPP_{n,k}}r(\lambda)K_{n,k}(\nu,\lambda).
    \end{align}
    Since $\PPP_{n,k}$ is finite, choose $\nu_0$ where $r$ attains its maximum. Then
    \[
    r(\nu_0) = \sum_{\lambda\in\PPP_{n,k}}r(\lambda) K_{n,k}(\nu_0,\lambda) \leq r(\nu_0) \sum_{\lambda\in\PPP_{n,k}}K_{n,k}(\nu_0,\lambda)= r(\nu_0).
    \]
    The equalities force $r(\lambda) = r(\nu_0)$ if $K_{n,k}(\nu_0,\lambda)>0.$ Repeating the same argument along paths of transitions with positive probability and using irreducibility, we obtain $r(\lambda) = r(\nu_0)$ for all $\lambda \in \PPP_{n,k}$. Since both $\theta$ and $\pi_{n,k}$ add up to $1$, $r(\nu_0)=1$ and consequently $\theta = \pi_{n,k}$.
\end{proof}
\begin{lemma}\label{lem:restorerule}
    Let $n\geq 2$ and  $1\leq k <n$. Suppose that $\lambda\in \PPP_{n,k}$ and $\nu\in \PPP_{n,k+1}$ satisfy $\lambda\RRR\nu$. Let $\widehat{\lambda}$ and $\widehat{\nu}$ have joint distribution $Q_{\lambda,\nu}$ from \autoref{cor:probdistributeonpairs}. If a chosen partition lies outside its prescribed class, replace it by the corresponding initial partition. Thus, we define the restoration formula
    \[
    \alpha=
    \begin{cases}
    \widehat\lambda,&\widehat\lambda_1=k,\\
    \lambda, &\widehat{\lambda}_1\neq k,
    \end{cases}
    \qquad
    \beta=
    \begin{cases}
    \widehat\nu,&\widehat\nu_1=k+1,\\
    \nu,&\widehat\nu_1\neq k+1.
    \end{cases}
    \]
    Then for every pair $(\widehat{\lambda},\widehat{\nu})$ satisfying $Q_{\lambda,\nu}(\widehat{\lambda},\widehat{\nu})>0$, the partitions $\alpha,\beta$ obtained by the restoration formula satisfy $\alpha\RRR\beta$.
\end{lemma}
\begin{proof}
    Fix a pair $(\widehat{\lambda},\widehat{\nu})$ with  $Q_{\lambda,\nu}(\widehat{\lambda},\widehat{\nu})>0$. By \autoref{cor:probdistributeonpairs}, we have $\widehat{\lambda}\RRR\widehat{\nu}$. The marginal identities for $Q_{\lambda,\nu}$ also imply $K_n(\lambda,\widehat{\lambda})>0$ and $K_n(\nu,\widehat{\nu})>0$. The $K_n$ transition removes one box and then adds one box, so
    \[
    |\widehat{\lambda}_1 - k | \leq 1,\qquad |\widehat{\nu}_1 - (k+1) | \leq 1.
    \]
    By construction, $\lambda_1' = k$ and $\nu_1' = k+1$. It therefore suffices to prove $\nu_i' \leq \lambda_i' +1$ for every $i\geq 2$. There are four cases.
    \begin{enumerate}
        \item Neither partition is restored. Then $(\alpha,\beta) = (\widehat{\lambda},\widehat{\nu})$, so the required inequality follows from $\widehat{\lambda}\RRR\widehat{\nu}$.
        \item Both partitions are restored. Then $(\alpha,\beta) = (\lambda,\nu)$, so the required inequality follows from $\lambda\RRR\nu$.
        \item Only $\widehat{\lambda}$ is restored. In this case $\widehat{\lambda}_1 \neq k$ and $\widehat{\nu}_1 = k+1$. Since $\widehat{\lambda}\RRR\widehat{\nu}$, $k+1 = \widehat{\nu}_1 \leq \widehat{\lambda}_1+1$. By the definition of $Q_{\lambda,\nu}$, we must have $\widehat{\lambda} = \lambda -e_i +e_{i'}$ for some $i,i'$ and the bound on the change in the first row forces $\widehat{\lambda}_1 = k+1$. Hence the move from $\lambda$ to $\widehat{\lambda}$ removes a box below the first row and adds a box to the first row. In particular, $\widehat{\lambda}_i \leq \lambda_i$ for all $i\geq 2$. Since $(\alpha,\beta) = (\lambda,\widehat{\nu})$, it follows that
        \[
        \beta_i = \widehat{\nu}_i \leq \widehat \lambda_i +1 \leq \lambda_i +1 = \alpha_i+1,\qquad i\geq 2.
        \]
        \item  Only $\widehat{\nu}$ is restored. By an argument similar to that in the previous case, we must have $\widehat{\nu}_1 = k$ and $\nu_i \leq \widehat{\nu}_i$ for all $i\geq 2$. Since $(\alpha,\beta) = (\widehat{\lambda},\nu)$, it follows that
        \[
        \beta_i = \nu_i \leq \widehat{\nu}_i \leq \widehat \lambda_i +1 = \alpha_i+1,\qquad i\geq 2.
        \]
    \end{enumerate}
\end{proof}
\begin{theorem}\label{thm:couplingRnk}
    For every $n\geq 2$ and $1 \leq k <n$, the conditional Plancherel distributions $\pi_{n,k}$ and $\pi_{n,k+1}$ admit a coupling $\theta$ supported on
    \[
    \RRR_{n,k} := \{ (\lambda,\nu)\in \PPP_{n,k}\times \PPP_{n,k+1} : \lambda\RRR\nu\}. 
    \]
\end{theorem}
\begin{proof}
   
    The set $\RRR_{n,k}$ is finite and nonempty since the pair $\lambda = (k, 1^{n-k}), \nu = (k+1, 1^{n-k-1})$ belongs to it.

    For each pair $(\lambda,\nu)\in \RRR_{n,k}$, choose $(\widehat{\lambda}, \widehat{\nu})$ with joint distribution $Q_{\lambda,\nu}$ from \autoref{cor:probdistributeonpairs}. Apply the rule in \autoref{lem:restorerule} and denote the resulting pair by $(\alpha,\beta)$. For every pair $(\widehat{\lambda}, \widehat{\nu})$ with $Q_{\lambda,\nu}(\widehat{\lambda}, \widehat{\nu}) > 0$, we have $(\alpha,\beta)\in\RRR_{n,k}$ by \autoref{lem:restorerule}. Let $\JJJ_{\lambda,\nu}:\RRR_{n,k}\rightarrow \RR_{\geq 0}$ be a row of the transition matrix $\mathcal{J}$ where $\JJJ_{\lambda,\nu}(\alpha,\beta)$ is the sum of $Q_{\lambda,\nu}(\widehat{\lambda},\widehat{\nu})$ over all pairs $(\widehat{\lambda},\widehat{\nu})$ whose restoration produces $(\alpha,\beta)$. Every chosen pair $(\widehat{\lambda},\widehat{\nu})$ produces exactly one restored pair in $\RRR_{n,k}$, so
    \[
    \sum_{(\alpha,\beta)\in\RRR_{n,k}}
    \JJJ_{\lambda,\nu}(\alpha,\beta)
    =
    \sum_{\widehat\lambda,\widehat\nu\in\PPP_n}
    Q_{\lambda,\nu}(\widehat\lambda,\widehat\nu)
    =1.
    \] Thus $\JJJ$ is a transition matrix on $\RRR_{n,k}$.
    
    Fix $\alpha\in \PPP_{n,k}$ with $\alpha\neq \lambda$. A move produces $\alpha$ as its first partition exactly when $\widehat{\lambda} = \alpha$. Hence
    \[
    \sum_{\beta:(\alpha,\beta)\in \RRR_{n,k}}\JJJ_{\lambda,\nu}(\alpha,\beta) = \sum_{\xi\in\PPP_n} Q_{\lambda,\nu}(\alpha,\xi)= K_{n}(\lambda,\alpha) = K_{n,k}(\lambda,\alpha),
    \]
    where the second equality follows from \autoref{cor:probdistributeonpairs} and the third follows from \autoref{def:restoremove}.

    If $\alpha = \lambda$,
    \begin{align*}
    \sum_{\beta:(\alpha,\beta)\in \RRR_{n,k}}\JJJ_{\lambda,\nu}(\alpha,\beta) 
    &= \sum_{\xi\in\PPP_n} Q_{\lambda,\nu}(\lambda,\xi)
    +\sum_{\zeta\in \PPP_n\setminus\PPP_{n,k}}\sum_{\xi\in\PPP_n} Q_{\lambda,\nu}(\zeta,\xi)\\
    &= K_n(\lambda,\lambda) + \sum_{\zeta\in \PPP_n\setminus \PPP_{n,k}}K_n(\lambda,\zeta)\\
    &= K_{n,k}(\lambda,\lambda).
    \end{align*}
    The same argument, summing over $\alpha\in \PPP_{n,k}$, gives $K_{n,k+1}(\nu,\beta)$ for all $\beta\in \PPP_{n,k+1}$. Thus, the transition matrix $\JJJ$ couples $K_{n,k}$ and $K_{n,k+1}$ as introduced in \autoref{subsec:coupledtrans}. 
    
    Since $\RRR_{n,k}$ is finite and nonempty, the transition matrix $\JJJ$ has a stationary probability distribution $\theta$. Define two distributions from $\theta$ by 
    \[
    \theta_1(\alpha) = \sum_{\beta:(\alpha,\beta)\in \RRR_{n,k}}\theta(\alpha,\beta),\quad \theta_2(\beta) = \sum_{\alpha:(\alpha,\beta)\in \RRR_{n,k}}\theta(\alpha,\beta).
    \]
    For each $\alpha\in \PPP_{n,k}$,
    \begin{align*}
    \theta_1(\alpha)
    &= \sum_{\beta:(\alpha,\beta)\in \RRR_{n,k}}\bigg( \sum_{(\lambda,\nu)\in \RRR_{n,k}}\theta(\lambda,\nu)\cdot\JJJ_{\lambda,\nu}(\alpha,\beta)\bigg)\\
    &= \sum_{(\lambda,\nu)\in \RRR_{n,k}}\bigg(\theta(\lambda,\nu)\cdot \sum_{\beta:(\alpha,\beta)\in\RRR_{n,k}}\JJJ_{\lambda,\nu}(\alpha,\beta)\bigg)\\
    &= \sum_{(\lambda,\nu)\in \RRR_{n,k}}\theta(\lambda,\nu)K_{n,k}(\lambda,\alpha)\\
    &= \sum_{\lambda\in \PPP_{n,k}}\theta_1(\lambda)K_{n,k}(\lambda,\alpha).
    \end{align*}
    Thus $\theta_1$ is stationary for $K_{n,k}$. By \autoref{lem:Kmk}, its unique stationary distribution is $\pi_{n,k}$, so $\theta_1 = \pi_{n,k}$. Similarly, we can derive $\theta_2 = \pi_{n,k+1}$.
    
    Finally, extending $\theta$ by zero outside $\RRR_{n,k}$ gives the required coupling on $\PPP_{n,k}\times \PPP_{n,k+1}$.
\end{proof}
\section{Proof of strict log-concavity}\label{sec:logconcavity}
Fix $n\ge1$. For $\lambda\in\PPP_n$, let $b(\lambda)$ and $d(\lambda)$ denote the probabilities that a single $K_n$ transition increases or decreases the first-row length by one, respectively. Thus,
\begin{align*}
    b(\lambda) &= \sum_{\eta\in\PPP_n:\eta_1= \lambda_1+1}K_{n}(\lambda,\eta),\\
    d(\lambda) &= \sum_{\eta\in\PPP_n:\eta_1= \lambda_1-1}K_{n}(\lambda,\eta).
\end{align*}
For $1\le k\le n$, define their expections under the conditional Plancherel measure $\pi_{n,k}$ by
\[
b_k = \mathbb{E}_{\pi_{n,k}}[b] = \sum_{\lambda\in \PPP_{n,k}}\pi_{n,k}(\lambda)b(\lambda),\quad d_{k} = \mathbb{E}_{\pi_{n,k}}[d] = \sum_{\lambda\in \PPP_{n,k}}\pi_{n,k}(\lambda)d(\lambda).
\]

First, we present an explicit formula of $b(\lambda)$ using the addition row probability $U_1(\lambda)$ defined in \autoref{eq:goup}.
\begin{lemma}\label{lem:strictineq}
    For every $\lambda \in \PPP_n$ with $n\geq 1$,
    \[
    b(\lambda) = \frac{(n+1)U_1(\lambda)-1}{n}.
    \]
\end{lemma}
\begin{proof}
    To obtain a partition $\eta$ with $\eta_1 = \lambda_1+1$ through a $K_n$ move, one must first remove a box below the first row and then add one to the first row. Thus, summing over all possible removals gives
    \[
    b(\lambda) = \sum_{\substack{\rho\nearrow \lambda\\ \rho_1 = \lambda_1}} \frac{f^\rho}{f^\lambda}\cdot \frac{f^{\rho+e_1}}{nf^\rho} = \sum_{\substack{\rho\nearrow \lambda\\ \rho_1 = \lambda_1}}\frac{f^{\rho+e_1}}{nf^\lambda}.
    \]
    By \autoref{eq:godown},
    \[
    f^{\lambda + e_1} = f^\lambda + \sum_{\substack{\rho\nearrow \lambda\\ \rho_1 = \lambda_1}}f^{\rho+e_1}.
    \]

    Therefore,
    \[
    b(\lambda) = \frac{f^{\lambda+e_1}-f^\lambda}{nf^\lambda} = \frac{(n+1)U_1(\lambda)-1}{n}.
    \]
\end{proof}

\begin{lemma}\label{lem:comparison}
    If $\lambda\in \PPP_{n,k}$, $\nu\in \PPP_{n,k+1}$ and $\lambda \RRR \nu$, then 
    \[
    b(\lambda) > b(\nu),\quad d(\lambda) \leq d(\nu).
    \]
\end{lemma}
\begin{proof}
    For every pair $(\alpha,\beta)\in\PPP_n\times\PPP_n$ with $Q_{\lambda,\nu}(\alpha,\beta)>0$, the marginal identities imply $K_n(\lambda,\alpha)>0$ and $K_n(\nu,\beta)>0$. Since a $K_n$ transition removes one box and then adds one box,
    \[
    k-1\leq\alpha_1\leq k+1,
    \qquad
    k\leq\beta_1\leq k+2.
    \]
    The support condition for $Q_{\lambda,\nu}$ also gives $\alpha\RRR\beta$, so $\beta_1\leq\alpha_1+1$.

    If $\beta_1=k+2$, then
    $k+1=\beta_1-1\leq\alpha_1\leq k+1$,
    and hence $\alpha_1=k+1$.
    Thus whenever the first row of $\nu$ grows, the first row of $\lambda$ also grows.
    Similarly, if $\alpha_1=k-1$, then
    $k\leq\beta_1\leq\alpha_1+1=k$,
    and hence $\beta_1=k$.
    Thus whenever the first row of $\lambda$ shrinks, the first row of $\nu$ also shrinks.

    Since the entries of $Q_{\lambda,\nu}$ are nonnegative, the marginal identities give
    \begin{align*}
    b(\nu)
    &=\sum_{\substack{\alpha,\beta\in\PPP_n\\ \beta_1=k+2}}Q_{\lambda,\nu}(\alpha,\beta)
    \leq\sum_{\substack{\alpha,\beta\in\PPP_n\\ \alpha_1=k+1}}Q_{\lambda,\nu}(\alpha,\beta)
    =b(\lambda),\\
    d(\lambda)
    &=\sum_{\substack{\alpha,\beta\in\PPP_n\\ \alpha_1=k-1}}Q_{\lambda,\nu}(\alpha,\beta)
    \leq\sum_{\substack{\alpha,\beta\in\PPP_n\\ \beta_1=k}}Q_{\lambda,\nu}(\alpha,\beta)
    =d(\nu).
    \end{align*}
    The strict inequality $b(\lambda) > b(\nu)$ follows from the explicit expression from \autoref{lem:strictineq} and the fact that $U_1(\lambda) > U_1(\nu)$ by \autoref{lem:DU} since row one is tight.
\end{proof}
\begin{proposition}\label{prop:monotone}
    For $1\leq k<n$,
    \[
    b_k > b_{k+1},\quad d_{k}\leq d_{k+1}.
    \]
    Moreover, $b_k > 0$ for $1 \leq k < n$ and $d_k > 0$ for $1 < k \leq n$.
\end{proposition}
\begin{proof}
    Fix $1\leq k < n$, and let $\theta(\lambda,\nu)$ be a coupling of $\pi_{n,k}$ and $\pi_{n,k+1}$ supported on $\RRR_{n,k}$ defined in \autoref{thm:couplingRnk}. For every $(\lambda,\nu)\in \RRR_{n,k}$, \autoref{lem:comparison} gives $b(\lambda) > b(\nu)$. Multiply by $\theta(\lambda,\nu)$ and sum over all pairs in $\RRR_{n,k}$:
    \[
    \sum_{(\lambda,\nu)\in \RRR_{n,k}}\theta(\lambda,\nu)b(\lambda)> \sum_{(\lambda,\nu)\in \RRR_{n,k}}\theta(\lambda,\nu)b(\nu).
    \]
    The left hand side equals \[\sum_{\lambda\in \PPP_{n,k}}\pi_{n,k}(\lambda)b(\lambda)= b_k,\] and the right hand side equals \[\sum_{\nu\in \PPP_{n,k+1}}\pi_{n,k+1}(\nu)b(\nu)= b_{k+1}.\]
    The same argument gives $d_k \leq d_{k+1}$.

    For positivity, now let $1\leq k \leq n$ and set $\lambda = (k ,1^{n-k})$. If $k<n$, it has a box below the first row. Removing its bottom box and adding a box to the first row gives $(k, 1^{n-k}) \rightarrow (k+1, 1^{n-k-1})$. This is a legal $K_n$ move with positive probability, so $b(\lambda) > 0$. Also $\pi_{n,k}(\lambda) >0$; hence, $b_k\geq \pi_{n,k}(\lambda)b(\lambda) > 0$. Similarly, when $k>1$, removing a box from the hook's first row and adding a new bottom box is a legal move in $K_n$. The same argument gives $d_k > 0$. 
\end{proof}

Now we relate the expectations with $a_{n,k}$ and $a_{n,k+1}$.
\begin{lemma}\label{lem:balance}
    For $1\leq k <n$,
    \[
    a_{n,k}b_k = a_{n,k+1}d_{k+1}.
    \]
\end{lemma}
\begin{proof}
    Take $\lambda\in \PPP_{n,k}$ and $\nu\in \PPP_{n,k+1}$. \autoref{eq:reversibility} gives 
    \[
    \pi_n(\lambda) K_{n}(\lambda,\nu) = \pi_n(\nu)K_{n}(\nu,\lambda).
    \]
    Since $\pi_n(\lambda) = \frac{(f^\lambda)^2}{n!}$, multiplying by $n!$ gives
    \[
    (f^\lambda)^2 K_{n}(\lambda,\nu)= (f^\nu)^2 K_{n}(\nu,\lambda).
    \]
    Summing over all such pairs gives
    \[
    \sum_{\lambda\in\PPP_{n,k}}\sum_{\nu\in \PPP_{n,k+1}}(f^\lambda)^2 K_n(\lambda,\nu)=
    \sum_{\nu\in\PPP_{n,k+1}}\sum_{\lambda\in \PPP_{n,k}}(f^\nu)^2 K_n(\nu,\lambda).
    \]
    A transition from $\lambda$ ends in $\PPP_{n,k+1}$ exactly when its first-row length increases by one. Thus the left hand side equals
    \[
    \sum_{\lambda\in\PPP_{n,k}}\sum_{\nu\in \PPP_{n,k+1}}(f^\lambda)^2 K_n(\lambda,\nu)=
    \sum_{\lambda\in\PPP_{n,k}}(f^\lambda)^2b(\lambda) = a_{n,k}\sum_{\lambda\in\PPP_{n,k}}\pi_{n,k}(\lambda)b(\lambda) = a_{n,k}b_k.
    \]
    A transition from $\nu$ ends in $\PPP_{n,k}$ exactly when its first-row length decreases by one. Thus the right hand side equals
    \[
    \sum_{\nu\in\PPP_{n,k+1}}\sum_{\lambda\in \PPP_{n,k}}(f^\nu)^2 K_n(\nu,\lambda)=
    \sum_{\nu\in\PPP_{n,k+1}}(f^\nu)^2d(\nu) = a_{n,k+1}\sum_{\nu\in\PPP_{n,k+1}}\pi_{n,k+1}(\nu)d(\nu) = a_{n,k+1}d_{k+1}.
    \]
\end{proof}
The final step uses the same principle as the criterion for log-concavity of stationary distributions of birth-and-death chains used by Whitt \cite[Theorem~3]{Whitt85}: nonincreasing birth rates and nondecreasing death rates imply log-concavity.

\begin{thm:Chen}
    For $1 < k <n$,
    \[
    a_{n,k}^2 > a_{n,k-1} a_{n,k+1}.
    \]
\end{thm:Chen}
\begin{proof}
    By \autoref{prop:monotone} and \autoref{lem:balance}, we have
    \[
    \frac{a_{n,k}}{a_{n,k-1}} = \frac{b_{k-1}}{d_k}> \frac{b_k}{d_k} \geq \frac{b_k}{d_{k+1}} = \frac{a_{n,k+1}}{a_{n,k}}.
    \]
     For $n = 1,2$, the sequences are $(1)$ and $(1,1)$ respectively, and there are no interior inequalities.
\end{proof}
\medskip

\noindent
\textbf{AI disclosure.}
The author used OpenAI's GPT-6 Astra for assistance with mathematical discussions and the revision of this manuscript. 
The author wrote the manuscript and takes full responsibility for the mathematical content and text.
\bibliographystyle{amsalpha}
	\bibliography{references}						
\end{document}